\documentclass[a4paper,11pt]{article}
\usepackage[margin=25mm]{geometry}
\usepackage{lmodern}
\usepackage{authblk}
\usepackage{xcolor}
\usepackage[T1]{fontenc}
\usepackage[utf8]{inputenc}
\usepackage{amsmath,amssymb,amsthm,mathtools}
\usepackage{graphicx,booktabs,array,microtype}
\usepackage{algorithm}
\usepackage{subcaption}
\usepackage{algpseudocode}
\usepackage{url}
\usepackage[colorlinks=true,
            linkcolor=blue!65!black,
            citecolor=blue!65!black,
            urlcolor=blue!65!black,
            bookmarksnumbered=true,
            pdfpagelabels=true]{hyperref}
\usepackage{bookmark}
\hypersetup{
  pdftitle={Data Informativity for Robust Feedback Rectification of Switched Linear Systems},
  pdfauthor={Philipp Schmitz, Hannes Gernandt, Maria C. Honecker, Karl Worthmann},
  pdfsubject={Single-column archive manuscript}
}
\allowdisplaybreaks[1]
\newtheorem{theorem}{Theorem}
\newtheorem{lemma}[theorem]{Lemma}
\newtheorem{corollary}[theorem]{Corollary}
\newtheorem{proposition}[theorem]{Proposition}
\theoremstyle{definition}
\newtheorem{definition}[theorem]{Definition}
\newtheorem{remark}[theorem]{Remark}
\newcommand{\R}{\mathbb R}
\newcommand{\C}{\mathbb C}
\newcommand{\Q}{\mathcal Q}
\newcommand{\norm}[1]{\lVert#1\rVert}
\newcommand{\im}{\operatorname{im}}
\newcommand{\rank}{\operatorname{rank}}
\newcommand{\diag}{\operatorname{diag}}

\newcommand{\spec}{\operatorname{spec}}

\newcommand{\bdir}{\beta_{\mathrm{dir}}}
\newcommand{\btri}{\beta}

\newcommand{\smat}[1]{\big[\begin{smallmatrix}#1\end{smallmatrix}\big]}
\title{Data Informativity for Robust Feedback Rectification of Switched Linear Systems\footnote{
\textbf{Funding.} HG was funded by
the Deutsche Forschungsgemeinschaft (DFG, German Research Foundation),
Project-ID 531152215 -- CRC 1701. PS  gratefully acknowledges funding by the Carl Zeiss
Foundation (VERNEDCT -- Project No. 2021-10-003).  KW is thankful for support
within the DFG research unit ALeSCo -- Active Learning in Systems and Control
(FOR 5785 -- Project-ID 535860958).
}
}

\author[1,2]{Hannes Gernandt\thanks{Corresponding author: \href{mailto:gernandt@uni-wuppertal.de}{gernandt@uni-wuppertal.de}.}}
\author[3]{Philipp Schmitz}
\author[4]{Maria C. Honecker}
\author[3]{Karl Worthmann}
\affil[1]{Port-Hamiltonian Institute, University of Wuppertal, Gaußstraße~20, 42119 Wuppertal, Germany}
\affil[2]{Fraunhofer Research Institution for Energy Infrastructures and Geotechnologies IEG, Gulbener Straße~23, 03046 Cottbus, Germany}
\affil[3]{Optimization-based Control Group, Technische Universit\"at Ilmenau, Weimarer Straße~25, 98693 Ilmenau, Germany}
\affil[4]{cigus GmbH, Magirus-Deutz-Straße~13, 89077 Ulm, Germany}
\date{}

\begin{document}
\maketitle

\begin{abstract}
    Feedback rectification is a method for switched linear systems that constructs a mode-dependent state feedback such that the closed-loop subsystems share a common eigenvector basis. 
    If the assigned eigenvalues lie inside the unit disk, this yields a common quadratic Lyapunov function and guarantees exponential stability under arbitrary switching. 
    In this paper, we study robust feedback rectification from input-state data affected by bounded disturbances. 
    Using the concept of data informativity, we derive necessary and sufficient conditions for approximate feedback rectification without requiring full row-rank of the collected input-state data. If the rectification error is smaller than the stability margin of the prescribed eigenvalues, we obtain global  
    exponential stability under arbitrary switching in the disturbance-free case, and practical global asymptotic stability under 
    bounded disturbances. 
    The approach is illustrated by a mass-spring-damper chain with varying stiffness, where prescribed poles are assigned in every mode from data that do not determine the system matrices uniquely. 
\end{abstract}
\section{Introduction}

Switched linear systems describe processes whose dynamics change between a finite number of operating modes, e.g., power converters~\cite{trenn2012switched}, aero-engine control~\cite{WangZS2016}, and precision motion systems~\cite{HeertjesN2012}. Since switching between stable subsystems may lead to instability, their stability analysis is more subtle than for a single linear system \cite{Liberzon.2003,ShortenEtAl2007,lin2009stability}. Stability can be recovered by restricting the switching signal, e.g., by average dwell-time conditions~\cite{Hespanha99}, or, under arbitrary switching, certified by a common quadratic Lyapunov function (CQLF). A CQLF is sufficient but not necessary; more general certificates use multiple or nonquadratic Lyapunov functions \cite{Branicky1998,ShortenEtAl2007}. 

Controller design for such systems has to reconcile two goals. \emph{Stability under arbitrary switching} is typically enforced by searching for feedback gains together with a CQLF, e.g., by LMI-based synthesis~\cite{Sun05} or $H_\infty$ control~\cite{Deaecto11}. In contrast, placing stable poles in every mode separately gives no guarantee under switching. The closed-loop matrices have different eigenvectors, and switching between them can cause large transient amplification or even instability. \emph{Performance of each mode}, on the other hand, is naturally specified by closed-loop poles. The CQLF-based methods above certify stability, but at best confine the closed-loop eigenvalues to a region, so the damping of each mode is bounded rather than prescribed. Both goals are partially addressed by constructive designs for single-input systems~\cite{WulWirSho09}, eigenvalue assignment with simultaneous triangularization~\cite{Hai16,SunShorten2005}, and predictive control~\cite{ZhangZB2016,Ju2026}.

Feedback rectification~\cite{HoneckerGWTR24,WulHR21} achieves both goals simultaneously. It builds on eigenstructure assignment: Moore's classical result~\cite{Moore76} characterizes the closed-loop eigenvectors compatible with prescribed distinct eigenvalues and thereby identifies the design freedom beyond pole placement; see \cite{Schmid14,padula2021eigenstructure} for extensions. Feedback rectification uses this freedom to choose mode-dependent gains such that all closed-loop matrices have the prescribed eigenvalues and a common eigenvector basis. In these coordinates, every closed-loop matrix is normal and has the prescribed poles as eigenvalues. Consequently, the norm in these coordinates defines a CQLF whose decay rate is set by the designer. No separate Lyapunov search is required, and this norm decreases monotonically without transient amplification. The same eigenstructure-based approach has been used for monotone tracking and non-overshooting output shaping~\cite{WulffHSR26}. Model-based rectification, however, requires the subsystem matrices and intersections of particular subspaces given by Moore's result~\cite{WulHR21,HoneckerGWTR24}.

Data-driven control avoids the need for a model and builds on  Willems' fundamental lemma~\cite{Willems2005} 
and on the informativity framework~\cite{van2020data}. 
When the data are disturbed, a single fitted model need not describe all systems consistent with the measurements. Robust designs therefore impose their requirements on the entire compatible model set, using tools from semidefinite programming such as the matrix $S$-lemma~\cite{vanWaarde2022}. 
Data-driven pole and eigenstructure assignment for a single linear system is studied in~\cite{celi23}. For switched systems, data-driven stabilization has been studied with online mode detection and dwell-time assumptions~\cite{eising2024data}, with a jointly designed switching signal~\cite{bianchi2025data}, and with output feedback~\cite{hu2023data}. In contrast, we assume the active mode to be known, do not design the switching signal, and require the feedback to realize prescribed closed-loop eigenvalues and common eigenvectors for every system compatible with the data. This extends the data-driven counterpart of exact feedback rectification from noise-free data~\cite{Schmitz2026} to disturbed data.

With disturbed data, exact rectification of all compatible systems is in general impossible, and the achievable accuracy has to be quantified. We measure it by a rectification error, the worst-case distance between the closed loop in the common eigenvector coordinates and the prescribed normal matrix, taken over all systems of a mode that are compatible with the data. This leads to the central trade-off of the paper. The prescribed poles fix the nominal decay rate, and the remaining stability margin is the budget available to absorb the model uncertainty.

We characterize when the data are informative for approximate feedback rectification. Using the $S$-lemma, we then derive a~semidefinite program for the worst-case rectification error and a sufficient LMI condition that computes the common eigenvectors and the gains jointly, avoiding the subspace  intersections from \cite{Schmitz2026}. The resulting CQLF certifies uniform exponential stability for all compatible systems under arbitrary switching and practical stability under bounded disturbances. For a lightly damped three-mass chain with switching stiffness, approximate feedback rectification from closed-loop data that do not determine the system matrices uniquely assigns the same poles in every mode and is certified for more than six times the noise level of a certainty-equivalent design, whereas mode-wise pole placement is destabilized by switching.

Section~\ref{sec:problem} introduces the system dynamics, data informativity, and approximate feedback rectification. Section~\ref{sec:data} characterizes approximate rectification from noisy 
data. Section~\ref{sec:design} develops 
semidefinite conditions, Section~\ref{sec:stability} derives the 
stability guarantees, and Section~\ref{sec:examples} presents a numerical example.

\noindent\textbf{Notation}.
$\mathbb{N}_0$ denotes the natural numbers including zero. 
The symbol $\norm{\,\cdot\,}$ denotes the Euclidean norm in $\mathbb{R}^n$
and the induced spectral matrix norm. The Moore--Penrose inverse is denoted by
$M^\dagger$. The set of eigenvalues of $M$ is $\spec M$. 
The image
and kernel are $\im M$ and $\ker M$. The inequalities $S\succeq0$ and
$S\succ0$ mean that $S\in\R^{n\times n}$ is symmetric positive semidefinite and
positive definite, respectively. For $S\succeq0$, $S^{1/2}$ denotes the unique symmetric positive semidefinite square root. The identity matrix in $\R^{n\times n}$ is denoted by $I_n$.

\section{Problem formulation}\label{sec:problem}
This section introduces the switched system dynamics subject to uncertainties and the problem of robust feedback rectification. 

\subsection{System, data, and disturbance model}
We first specify the measured input-state transitions and the set of system matrices that are consistent with them. 
For an arbitrary switching signal $\sigma:\mathbb N_0\to\Q:=\{1,\ldots,Q\}$ which is not known a priori, we consider the discrete-time system 
\begin{align}
\label{eq:switched_system}
    x_{k+1}=A_{\sigma(k)}x_k+B_{\sigma(k)}u_k,
\end{align}
with unknown matrices $A_q\in\R^{n\times n}$, and $B_q\in\R^{n\times m}$, where $x_k\in\R^n$, $u_k\in\R^m$ denote the state and the control at time~$k$ and $q = \sigma(k)$ the mode.

For each mode $q\in\Q$, collect $T_q$ state transitions in the 
data matrices
$X_q\in\R^{n\times T_q}$, $U_q\in\R^{m\times T_q}$, and
$X_q^+\in\R^{n\times T_q}$. Their columns contain the current state and 
input as well as the successor state. 
The measured state transitions satisfy 
\begin{equation}
\label{eq:system_dynamics}
     X_q^+ = A_qX_q+B_qU_q+W_q,\quad 
    W_q\in\R^{n\times T_q},\quad  
\end{equation}
with $\norm{W_q}\le\delta_q$,
where $\delta_q \ge 0$ is a known upper bound on the disturbance during data collection. 

The set of all \emph{compatible mode-$q$ systems with the measurements} $(X_q, U_q, X_q^+)$ is 
\[
 \Sigma_q^{\delta_q}:=\{(A,B)\in\R^{n\times n}\times\R^{n\times m}:\norm{X_q^+-AX_q-BU_q}\le\delta_q\}.
\]
This set contains all matrix pairs $(A,B)$ that could describe
the mode-$q$ dynamics consistently with the available measurements
and the prescribed error bound $\delta_q$.
In particular, since $\delta_q$ bounds the measurement error,
the inclusion $(A_q,B_q)\in\Sigma_q^{\delta_q}$ holds for the true system matrices, so $\Sigma_q^{\delta_q}\neq\emptyset$.
For the stacked data matrix
$Z_q:=\smat{X_q\\U_q}
\in\R^{(n+m)\times T_q}$ with $r_q:=\rank Z_q$, we choose $K_q\in\R^{T_q\times r_q}$ with orthonormal columns spanning $\im Z_q^\top$ and define the orthogonal projector
\[
 \Pi_q:=
 Z_q^{\dagger}Z_q
 =K_qK_q^\top
 \in\R^{T_q\times T_q}.
\]
Then $K_q^\top K_q=I_{r_q}$ and $(I_{T_q}-\Pi_q)K_q=0$.
In practice, $\delta_q$ is an assumption on the disturbance. The projector $\Pi_q$ allows us to check this assumption against the data.
\begin{lemma}
\label{lem:nonempty}
    $\Sigma_q^{\delta_q}\neq\emptyset$ if and only if
    \begin{equation}
    \label{eq:bound_proj}
        \Vert X_q^+(I_{T_q}-\Pi_q)\Vert\le\delta_q.
    \end{equation}
    In the case $\delta_q =0$ this is equivalent to
    \begin{equation*}
        \rank\begin{bmatrix}
            X_q\\U_q\\ X_q^+
        \end{bmatrix} = \rank\begin{bmatrix}
            X_q\\U_q
        \end{bmatrix}.
    \end{equation*}
\end{lemma}
\begin{proof}
Let $(A,B)\in\Sigma_q^{\delta_q}$, i.e.\ $\Vert W_q\Vert\leq \delta_q$ for $W_q := X_q^+ - \begin{bmatrix}A&\!\!\!\!B\end{bmatrix} Z_q$. Since $Z_q = Z_q\Pi_q$, we have
    $\begin{bmatrix}A&\!\!\!\!B\end{bmatrix} Z_q (I_{T_q} - \Pi_q)=0$ and
    \begin{equation*}
        X_q^+ (I_{T_q}-\Pi_q) \!=\! \big(W_q + \begin{bmatrix}
            A &\!\!\!\! B
        \end{bmatrix} Z_q\big) (I_{T_q} -\Pi_q)\!=\! W_q(I_{T_q}-\Pi_q).
    \end{equation*}
    As $I_{T_q}-\Pi_q$ is an orthogonal projection, $\Vert I_{T_q}-\Pi_q\Vert \leq 1$ and
    \[\Vert X_q^+ (I_{T_q}-\Pi_q)\Vert \leq \Vert W_q\Vert \leq \delta_q.\]

    Conversely, given \eqref{eq:bound_proj}, $\begin{bmatrix}A&\!\!\!\!B\end{bmatrix}:=X_q^+ Z_q^\dagger$ is an element of $\Sigma_{q}^{\delta_q}$. The statement for $\delta_q=0$ was shown in \cite{Schmitz2026}.
\end{proof}

For noise-free data, $(A_q,B_q)$ is uniquely identifiable if and only if $\rank Z_q=n+m$, which holds, e.g., if the input is persistently exciting of order $n+1$~\cite{Willems2005}. For disturbed data, full row rank makes the set $\Sigma_q^{\delta_q}$ bounded. Both are not required here.

\subsection{Prescribed closed-loop form and common contraction}
We next state the desired closed-loop eigenstructure of \eqref{eq:switched_system} and the separate condition that will be used to prove stability under arbitrary switching. While the switching signal~$\sigma$ is arbitrary and 
not known a priori, the active mode~$q$ can be tracked instantaneously for applying the mode-dependent feedback law $u_k=F_{\sigma(k)}x_k$ with $F_q\in\R^{m\times n}$ for $q = \sigma(k)$. 
For each mode $q\in\mathcal{Q}$, let the desired closed-loop eigenvalues $\lambda_{q,1},\ldots,\lambda_{q,n}\in\C$ be given that are closed under complex conjugation, counting multiplicities.

Feedback rectification of \eqref{eq:switched_system} seeks state feedbacks $F_q$, $q\in\mathcal{Q}$, that place these eigenvalues such that all closed-loop matrices $A_q+B_qF_q$ share a common set of eigenvectors $\{v_i\}_{i=1}^n$.
More precisely, we require a nonsingular $V\in\R^{n\times n}$ such that
\begin{align}
\nonumber
&V^{-1}(A_q+B_qF_q)V=\Lambda_q=\diag(\Lambda_q^{(j)})_j \quad \text{for all $q\in\mathcal{Q}$, with}\\
& \Lambda_q^{(j)}=\begin{cases}
    \lambda_{q,j} & \text{if $\lambda_{q,j}\in\mathbb{R}$,}\\
    \smat{a_{q,j}&b_{q,j}\\-b_{q,j}&a_{q,j}} & \text{if $\lambda_{q,j}=a_{q,j}+\mathrm{i}b_{q,j}$, $b_{q,j}>0$,}
\end{cases}
\label{def_Lambda_q}
\end{align}
Each $2\times 2$ block represents a complex-conjugate pair $a_{q,j}\pm \mathrm i b_{q,j}$. These eigenvalues occur in conjugate pairs since all matrices are real. 
Note that with this structure, the $\Lambda_q$ are normal, and we assume
$\norm{\Lambda_q}<1$ for all $q\in\Q$. By normality, this is equivalent to
all prescribed eigenvalues $\lambda_{q,i}$ lying in the open unit disk. For real diagonal targets, the
columns of $V$ are common eigenvectors. For a complex pair, the corresponding
two columns span a common real invariant subspace.

\begin{definition}[Approximate feedback rectification]\label{def:rectification}
For every $q\in\Q$, let $\delta_q\ge0$, prescribed closed-loop poles $\lambda_{q,1},\ldots,\lambda_{q,n}\in\C$ with $\Lambda_q$ given by \eqref{def_Lambda_q}, and a tolerance $\varepsilon_q\ge0$ be given. The data $(X_q,U_q,X_q^+)_{q\in\Q}$ are called \emph{informative for approximate $(\varepsilon_q)_{q\in\Q}$-feedback rectification with poles $(\Lambda_q)_{q\in\Q}$} if there exist one nonsingular $V\in\R^{n\times n}$, common to all modes, and, for every $q\in\Q$, a gain $F_q\in\R^{m\times n}$ such that 
\[
 \sup_{(A,B)\in\Sigma_q^{\delta_q}}
 \norm{V^{-1}(A+BF_q)V-\Lambda_q}\le\varepsilon_q.
\]
If $\varepsilon_q=0$ for all $q\in\Q$, then we say that the data are informative for \emph{exact} feedback rectification.
\end{definition}
Here, the bound $\delta_q$ describes the disturbance in the data, and the tolerance $\varepsilon_q$ the accuracy required of the design. For noise-free data, i.e., $\delta_q=0$, and $\varepsilon_q=0$, Definition~\ref{def:rectification} reduces to the notion in~\cite{Schmitz2026}. With $\varepsilon_q>0$, also systems that are not exactly rectifiable can be handled. For disturbed data, exact rectification is in general impossible, so the case of interest is $\delta_q>0$ and $\varepsilon_q>0$. Then $\varepsilon_q$ can be interpreted as the part of the stability margin $1-\|\Lambda_q\|$ of the prescribed poles that is consumed by the uncertainty in the data.

\section{Robust data-driven feedback rectification}\label{sec:data}

This section separates the disturbance fixed by the measurements from the part that remains free and then characterizes exactly which feedback gains can satisfy
a finite robust error bound. 
For given data $(X_q,U_q,X_q^+)$ the matrix
\begin{align*}
 \Omega_q:=\delta_q^2I_n-X_q^+(I_{T_q}-\Pi_q)(X_q^+)^\top\in\R^{n\times n}
\end{align*}
measures how much of the disturbance bound remains after the part $X_q^+(I_{T_q}-\Pi_q)$ of the disturbance has been taken into account.

\begin{lemma}[Disturbance parametrization]\label{lem:reduction}
We have $\Sigma_q^{\delta_q}\neq\emptyset$ if and only if  
$\Omega_q\succeq0$. In this case, $W_q\in \R^{n\times T_q}$ satisfies $X_q^+ = AX_q + B U_q + W_q$ for some $(A,B)\in\Sigma_q^{\delta_q}$ if and only if  there is $\Delta_q\in\R^{n\times r_q}$ with $\norm{\Delta_q}\le1$ such that 
\begin{equation}
\label{eq:Wq}
    W_q=X_q^+(I_{T_q}-\Pi_q)+\Omega_q^{1/2}\Delta_qK_q^\top.
\end{equation}
\end{lemma}
\begin{proof}
Since $I_{T_q}-\Pi_q$ is symmetric and idempotent, $\Omega_q\succeq 0$ is equivalent to \eqref{eq:bound_proj} and hence to $\Sigma_q^{\delta_q}\neq\emptyset$ by Lemma~\ref{lem:nonempty}. Assume $W_q = X_q^+ - \begin{bmatrix}A&B\end{bmatrix} Z_q$ for some $(A,B)\in\Sigma_q^{\delta_q}$, that is $W_qW_q^\top\preceq \delta_q^2 I_{n}$. With $Z_q\Pi_q=Z_q$, one sees $W_q (I_{T_q}-\Pi_q) = X_q^+ (I_{T_q}-\Pi_q)$. Therefore,
\begin{equation}
\label{eq:WW}
    W_q\Pi_q W_q^\top = W_qW_q^\top - X_q^+ (I_{T_q}-\Pi_q) (X_q^+)^\top \preceq \Omega_q.
\end{equation}
Since both matrices $W_q\Pi_q W_q^\top$ and $\Omega_q$ are symmetric positive semidefinite, this yields $\ker \Omega_q\subset \ker (W_q\Pi_q W_q^\top)$ and $\im (W_q\Pi_q) = \im (W_q\Pi_q W_q^\top) \subset \im \Omega_q = \im \Omega_q^{1/2}$. Set $\Delta_q := (\Omega_q^{1/2})^\dagger W_q K_q$. Then $\Omega_q^{1/2} \Delta_q K_q^\top = \Omega_q^{1/2}(\Omega_q^{1/2})^\dagger W_q\Pi_q = W_q\Pi_q$. Together with $W_q=W_q(I_{T_q}-\Pi_q)+W_q\Pi_q$ this shows \eqref{eq:Wq}. Since $\Omega_q^{1/2}$ is symmetric, $(\Omega_q^{1/2})^\dagger\Omega_q^{1/2}=\Omega_q^{1/2}(\Omega_q^{1/2})^\dagger$ is an orthogonal projection. With~\eqref{eq:WW} we find
\begin{equation*}
   \Delta_q \Delta_q^\top  = (\Omega_q^{1/2})^\dagger W_qK_qK_q^\top W_q^\top (\Omega_q^{1/2})^\dagger\preceq  (\Omega_q^{1/2})^\dagger \Omega_q (\Omega_q^{1/2})^\dagger = (\Omega_q^{1/2})^\dagger\Omega_q^{1/2} \preceq I_n.
\end{equation*}
This shows $\Vert \Delta_q \Vert \leq 1$.

Conversely, assume \eqref{eq:Wq} with $\Vert \Delta_q \Vert\leq 1$. Define \[\begin{bmatrix}A&\!\!\!B\end{bmatrix}:=(X_q^+ -W_q)Z_q^\dagger. 
\]
Then $\begin{bmatrix}A&\!\!\!B\end{bmatrix}Z_q = (X_q^+ -W_q)\Pi_q$. Thus,
$X_q^+ - \begin{bmatrix}A&\!\!\!B\end{bmatrix} Z_q = X_q^+ (I_{T_q}-\Pi_q) + \Omega_q^{1/2} \Delta_q K_q^\top = W_q$. Moreover, using $(I_{T_q}-\Pi_q)K_q=0$ and $K_q^\top K_q=I_{r_q}$, $W_q W_q^\top = X_q^+ (I_{T_q}-\Pi_q) (X_q^+)^\top + \Omega_q^{1/2} \Delta_q \Delta_q^\top \Omega_q^{1/2}$. Since $\norm{\Delta_q}\le1$, we have $\Delta_q \Delta_q^\top\preceq I_n$ and therefore $\Omega_q^{1/2} \Delta_q \Delta_q^\top \Omega_q^{1/2}\preceq \Omega_q$. Consequently, $W_qW_q^\top \preceq X_q^+ (I_{T_q}-\Pi_q) (X_q^+)^\top + \Omega_q = \delta_q^2 I_n$ implying $(A,B)\in\Sigma_q^{\delta_q}$.
\end{proof}
In particular, if $\delta_q=0$,
consistency forces $X_q^+(I_{T_q}-\Pi_q)=0$ and hence $\Omega_q=0$. Lemma~\ref{lem:reduction}
then reduces to the noise-free case $W_q=0$.

In the following, we present a first characterization of data informativity for robust feedback rectification. 
\begin{proposition}[Data informativity]
\label{thm:informativity}
The data are informative for
approximate $(\varepsilon_q)_{q\in\Q}$-feedback rectification with poles $(\Lambda_q)_{q\in\Q}$ if and only if
there exist a nonsingular $V\in\R^{n\times n}$ and matrices
$H_q\in\R^{r_q\times n}$ such that
\begin{equation}\label{eq:Hcondition}
 X_q K_qH_q=V,
\end{equation}
and
\begin{equation}
\label{eq:worst_case_error}
 \max_{\substack{\Delta_q\in\R^{n\times r_q}\\\norm{\Delta_q}\le1}}
 \norm{V^{-1}(X_q^+K_q-\Omega_q^{1/2}\Delta_q)H_q-\Lambda_q}
 \le\varepsilon_q
\end{equation}
for all $q\in\Q$. The corresponding gain is $F_q=U_qK_qH_qV^{-1}\in\R^{m\times n}$.
\end{proposition}

\begin{proof}
Suppose that the data are informative for
approximate $(\varepsilon_q)_{q\in\Q}$-feedback rectification with poles $(\Lambda_q)_{q\in\Q}$. Fix $(A^0_q,B^0_q)\in\Sigma_{q}^{\delta_q}$ and let
$L=\begin{bmatrix}L_A&L_B\end{bmatrix}\in\R^{n\times(n+m)}$ with $LZ_q=0$. For $t\in\R$ let $(A_q^t, B_q^t):=(A_q^0 + t L_A, B_q^0 + t L_B)$. Then $X_q^+ - \begin{bmatrix}A_q^t&B_q^t\end{bmatrix} Z_q = X_q^+ - \begin{bmatrix}A_q^0&B_q^0\end{bmatrix} Z_q$ and hence $(A_q^t, B_q^t)\in\Sigma_q^{\delta_q}$ for all $t\in\R$. Moreover,
\begin{equation*}
    \varepsilon_q\geq{} \Vert V^{-1} (A_q^t + B_q^t F_q) V - \Lambda_q\Vert \geq{}t\Vert V^{-1} L\smat{I_n\\F_q} V\Vert
    - \Vert V^{-1} (A_q^0 + B_q^0 F_q) V - \Lambda_q\Vert
\end{equation*}
for all $t\geq 0$, which implies $L\smat{I_n\\F_q} = 0$. This shows $\ker Z_q^\top \subset \ker\begin{bmatrix}I_n&F_q^\top\end{bmatrix}$ and equivalently $\im \smat{I_n\\F_q} \subset \im Z_q=\im Z_qK_q$. Therefore, there exists $H_q\in \R^{r_q\times n}$ satisfying
\begin{equation}
\label{eq:G}
    Z_q K_qH_q = \begin{bmatrix}I_n\\F_q\end{bmatrix} V.
\end{equation}
In particular, $X_qK_qH_q=V$ and $F_q=U_qK_qH_qV^{-1}$.

Let $\Delta_q\in\R^{n\times r_q}$ with $\norm{\Delta_q}\le1$. By
Lemma~\ref{lem:reduction} there is $(A_q,B_q)\in\Sigma_q^{\delta_q}$ with
$W_q:=X_q^+-\begin{bmatrix}A_q&B_q\end{bmatrix}Z_q=X_q^+(I_{T_q}-\Pi_q)+\Omega_q^{1/2}\Delta_qK_q^\top$.
With \eqref{eq:G}, $(I_{T_q}-\Pi_q)K_q=0$ and $K_q^\top K_q=I_{r_q}$ we obtain
\begin{equation}
\label{eq:closedloop}
(A_q+B_qF_q)V=\begin{bmatrix}A_q&B_q\end{bmatrix}Z_qK_qH_q=(X_q^+-W_q)K_qH_q=(X_q^+K_q-\Omega_q^{1/2}\Delta_q)H_q.
\end{equation}
Hence
$\norm{V^{-1}(X_q^+K_q-\Omega_q^{1/2}\Delta_q)H_q-\Lambda_q}
=\norm{V^{-1}(A_q+B_qF_q)V-\Lambda_q}\le\varepsilon_q$.

Conversely, let $V$ and $H_q$ satisfy the conditions of the proposition and
$F_q=U_qK_qH_qV^{-1}$. Then \eqref{eq:G} holds. Let
$(A_q,B_q)\in\Sigma_q^{\delta_q}$. By Lemma~\ref{lem:reduction} there is
$\Delta_q\in\R^{n\times r_q}$ with $\norm{\Delta_q}\le1$ such that
$W_q:=X_q^+-\begin{bmatrix}A_q&B_q\end{bmatrix}Z_q$ satisfies \eqref{eq:Wq}. Thus
\eqref{eq:closedloop} holds and
$\norm{V^{-1}(A_q+B_qF_q)V-\Lambda_q}\le\varepsilon_q$.
\end{proof}

\noindent By the first part of the proof of Proposition~\ref{thm:informativity}, every gain with finite worst-case error is of the form
\begin{equation*}
F_q=U_qK_qH_qV^{-1},\qquad H_q\in\mathbb R^{r_q\times n},\ X_qK_qH_q=V.
\end{equation*}
Their transformed closed-loop matrices form the set
\begin{equation*}
 \{V^{-1}X_q^+K_qH_q-V^{-1}\Omega_q^{1/2}\Delta H_q:
        \Delta\in\mathbb R^{n\times r_q},\ \|\Delta\|\leq1\}.
\end{equation*}

Since $V$ is nonsingular, \eqref{eq:Hcondition} implies $\rank X_q=n$. Hence, the assumption $\rank X_q=n$ in Theorem~\ref{thm:lmi-rectification} below is necessary for informativity.

\section{An LMI characterization of robust feedback rectifiability}\label{sec:design}
This section provides a norm test for the worst-case rectification error and derives an LMI certificate that computes the common eigenvectors and the feedback gains.

For $a,t\geq0$, $N,V\in\mathbb R^{n\times n}$ with $V$ nonsingular, and $H\in\mathbb R^{r_q\times n}$, define
\begin{equation}\label{eq:LMI}
 \mathcal M_q(a,N,V,H,t)\!:=\!
 \begin{bmatrix}
 aI_n\!&N^\top&\!H^\top\\
 N\!&I_n-tV^{-1}\Omega_qV^{-\top}&\!0\\
 H\!&0&\!tI_{r_q}
 \end{bmatrix}\!, 
\end{equation}
which is a symmetric matrix of size $(2n+r_q)\times(2n+r_q)$.
\begin{lemma}\label{lem:norm}
Let $H\neq0$, $V$ be invertible, and $\Omega_q\neq0$. Then
\begin{align}
    \label{eq:max_estimate}
 \max_{\substack{\Delta\in\mathbb R^{n\times r_q}\\\|\Delta\|\leq1}}
       \|N-V^{-1}\Omega_q^{1/2}\Delta H\|^2\leq a
\end{align}
is equivalent to $\mathcal M_q(a,N,V,H,t)\succeq0$ for some $t>0$. 
\end{lemma}

\begin{proof}
We abbreviate $C_q:=V^{-1}\Omega_q^{1/2}$. 
For any $x\in\mathbb R^n$ and
$\Delta\in\mathbb R^{n\times r_q}$ with $\|\Delta\|\leq1$,
the vector $y=\Delta Hx$ satisfies $\|y\|\leq\|Hx\|$.
Conversely, every $y\in\mathbb R^n$ satisfying this
inequality can be represented as $y=\Delta Hx$ with
$\|\Delta\|\leq1$. Indeed, if $Hx\neq0$, we may choose
$\Delta=\frac{y(Hx)^\top}{\|Hx\|^2}$, whereas $Hx=0$ implies $y=0$.

Consequently, the norm inequality \eqref{eq:max_estimate} is equivalent to
\begin{equation}\label{eq:implication}
 \|Hx\|^2-\|y\|^2\geq0
 \,\Longrightarrow\,
 a\|x\|^2-\|Nx-C_qy\|^2\geq0
\end{equation}
for all $x,y\in\mathbb R^n$.

We next apply the scalar $S$-lemma to \eqref{eq:implication}, which is possible since 
$H\neq0$ and, hence, choosing $x$ with $Hx\neq0$ and $y=0$
makes the condition on the left in \eqref{eq:implication} 
strictly positive.
Therefore, by the $S$-lemma
\cite[Section III-A]{vanWaarde2022},
\eqref{eq:implication} holds if and only if
there exists $\tau\geq0$ such that
\[
 a\|x\|^2-\|Nx-C_qy\|^2
 -\tau(\|Hx\|^2-\|y\|^2)\geq0
\]
for all $x,y\in\mathbb R^n$.
Equivalently,
\begin{equation}\label{eq:Slemma}
 \begin{bmatrix}
 aI_n-\tau H^\top H&0\\
 0&\tau I_n
 \end{bmatrix}
 -
 \begin{bmatrix}N^\top\\-C_q^\top\end{bmatrix}
 \begin{bmatrix}N&-C_q\end{bmatrix}
 \succeq0.
\end{equation}
The lower-right block of this inequality gives
$\tau I_n-C_q^\top C_q\succeq0$.
Since $C_q\neq0$, we necessarily have $\tau>0$.

It remains to transform \eqref{eq:Slemma} into $\mathcal{M}_q\succeq0$. Taking the Schur complement of the
middle block $I_n$ shows that \eqref{eq:Slemma}
is equivalent to 
\[
 \begin{bmatrix}
 aI_n-\tau H^\top H&N^\top&0\\
 N&I_n&-C_q\\
 0&-C_q^\top&\tau I_n
 \end{bmatrix}\succeq0.
\]
Moreover, a repeated application of Schur complements to the last and first block and using $t=\tau^{-1}>0$ leads to $\mathcal M_q(a,N,V,H,t)\succeq0$.
Since all transformations are equivalences, this proves the lemma. 
\end{proof}

The following theorem is our main result and presents a certificate for approximate feedback rectification. 
\begin{theorem}[Rectifiability certificate]
\label{thm:lmi-rectification} 
Assume $\rank X_q=n$ for every $q\in\Q$. Let $r_q$ and $K_q$ be as in Section~\ref{sec:problem}. Let $L_q\in\R^{r_q\times(r_q-n)}$ have columns spanning $\ker(X_qK_q)$. For $V\in\R^{n\times n}$ and $Y_q\in\R^{(r_q-n)\times n}$, set \[ H_q=(X_qK_q)^\dagger V+L_qY_q,\qquad R_q=X_q^+K_qH_q-V\Lambda_q. \] For prescribed tolerances $\varepsilon_q\geq0$, the data are informative for approximate $(\varepsilon_q)_{q\in\Q}$-feedback rectification with poles $(\Lambda_q)_{q\in\Q}$ if and only if there exist a nonsingular $V\in\R^{n\times n}$, matrices $Y_q$, and, whenever $\Omega_q\neq0$, scalars $t_q>0$, such that \begin{equation}\label{eq:unrestricted-rectification-test} 
\mathcal M_q(
\varepsilon_q^2,V^{-1}R_q,V,H_q,t_q)
\succeq0,
\end{equation} 
for each $q$ with $\Omega_q\neq0$. For $\Omega_q=0$, replace \eqref{eq:unrestricted-rectification-test} by
\begin{equation}\label{eq:unrestricted-rectification-zero} 
\begin{bmatrix} \varepsilon_q^2I_n&R_q^\top\\ R_q&VV^\top \end{bmatrix}\succeq0. 
\end{equation} 
The corresponding gains are $F_q=U_qK_qH_qV^{-1}\in\R^{m\times n}$. \end{theorem} \begin{proof} Since $X_q=X_qK_qK_q^\top$ and $\rank X_q=n$, $X_qK_q$ has full row-rank. Hence $H_q$ parametrizes precisely all solutions of $X_qK_qH_q=V$. By Proposition~\ref{thm:informativity}, data informativity is equivalent to the existence of nonsingular $V$ and such matrices $H_q$ satisfying \[ \max_{\|\Delta\|\le1} \|V^{-1}(R_q-\Omega_q^{1/2}\Delta H_q)\| \le\varepsilon_q, \] where $\Delta\in\R^{n\times r_q}$. Since $X_qK_qH_q=V$ is nonsingular, $H_q\neq0$. For $\Omega_q\neq0$, apply Lemma~\ref{lem:norm} with $N=V^{-1}R_q$ and $a=\varepsilon_q^2$. This gives exactly \eqref{eq:unrestricted-rectification-test}. For $\Omega_q=0$, the robust inequality reduces to $\|V^{-1}R_q\|\le\varepsilon_q$, which is equivalent by a Schur complement to \eqref{eq:unrestricted-rectification-zero}. All transformations are equivalences, proving both directions. 
\end{proof}

Since $V$ enters \eqref{eq:unrestricted-rectification-test} through $V^{-1}$, Theorem~\ref{thm:lmi-rectification} is not convex in $V$. Fixing a reference basis $V_0$ yields a sufficient condition that is an LMI in all remaining variables.

\begin{corollary}[LMI certificate]
\label{cor:lmi-reference} 
Under the assumptions of Theorem~\ref{thm:lmi-rectification} with $H_q$ and $R_q$, fix any nonsingular $V_0\in\R^{n\times n}$ and set $\widetilde R_q=V_0^{-1}R_q$.
If there exist $V$, $Y_q$, and $t_q\geq0$ such that 
\begin{equation}\label{eq:lmi-coordinate-condition} V_0^{-1}V+(V_0^{-1}V)^\top\succeq2I_n 
\end{equation} 
and, for every $q$ with $\Omega_q\neq0$, 
\begin{equation*}
\mathcal M_q(
\varepsilon_q^2,\widetilde R_q,V_0,H_q,t_q)
\succeq0
\end{equation*} 
and, for every $q$ with $\Omega_q=0$, $\smat{\varepsilon_q^2I_n&\widetilde R_q^\top\\ \widetilde R_q&I_n}\succeq0$.
Then the data are informative for approximate $(\varepsilon_q)_{q\in\Q}$-feedback rectification with poles $(\Lambda_q)_{q\in\Q}$, and $F_q=U_qK_qH_qV^{-1}$. 
\end{corollary} 
\begin{proof} Condition~\eqref{eq:lmi-coordinate-condition} gives 
\begin{align*} 
V_0^{-1}V(V_0^{-1}V)^\top-I_n=(V_0^{-1}V\!-\!I_n)(V_0^{-1}V\!-\!I_n)^\top \!+\!(V_0^{-1}V+(V_0^{-1}V)^\top\!-\!2I_n)\succeq0. 
\end{align*}
Hence $V_0^{-1}V$ and $V=V_0(V_0^{-1}V)$ are nonsingular and $VV^\top\succeq V_0V_0^\top$. 
Let $\Omega_q\neq0$. A congruence transformation with $\operatorname{diag}(I_n,V_0,I_{r_q})$ shows that $\mathcal M_q(\varepsilon_q^2,\widetilde R_q,V_0,H_q,t_q)\succeq0$ is equivalent to
\[
\begin{bmatrix} \varepsilon_q^2I_n&R_q^\top&H_q^\top\\ R_q&V_0V_0^\top-t_q\Omega_q&0\\ H_q&0&t_qI_{r_q} \end{bmatrix}\succeq0.
\]
Replacing the middle diagonal block by the larger matrix $VV^\top-t_q\Omega_q$ preserves positive semidefiniteness. A congruence transformation with $\operatorname{diag}(I_n,V^{-1},I_{r_q})$ then yields
\[
\mathcal M_q(\varepsilon_q^2,V^{-1}R_q,V,H_q,t_q)\succeq0,
\]
that is, \eqref{eq:unrestricted-rectification-test} holds, where $t_q>0$ since $X_qK_qH_q=V$ gives $H_q\neq0$ and $t_q=0$ would force $H_q=0$. For $\Omega_q=0$, a congruence transformation with $\operatorname{diag}(I_n,V_0)$ together with $V_0V_0^\top\preceq VV^\top$ yields \eqref{eq:unrestricted-rectification-zero}. The assertion follows from Theorem~\ref{thm:lmi-rectification}.
\end{proof}

\begin{remark}\label{rem:exact}
For $\delta_q=\varepsilon_q=0$,
Proposition~\ref{thm:informativity} recovers the
characterization of~\cite{Schmitz2026}, where $V$
is constructed using intersections of the subspaces
$X_q\ker(X_q^+ \Pi_q-\lambda_{q,i}X_q)$. Here, $X_q^+ \Pi_q$ is used instead of $X_q^+$ as in~\cite{Schmitz2026} to account for the disturbance $W_q$. 
Alternatively, $V$ and the feedback gains can be
computed using Corollary~\ref{cor:lmi-reference},
provided its conditions are feasible for the
chosen reference matrix $V_0$. For fixed nonsingular $V_0$, these conditions
are linear matrix inequalities in $V$, $Y_q$,
$t_q$, and $\varepsilon_q^2$.
\end{remark}

\section{Robust stability under arbitrary switching}\label{sec:stability}
We now show that a robust rectification based on collected data can be used to prove robust decay bounds for system operation under bounded disturbances
\begin{align}
\label{eq:with_noise}
    x_{k+1}=(A_{\sigma(k)}+B_{\sigma(k)}F_{\sigma(k)})x_k+d_k.
\end{align}

\begin{proposition}[Robust decay]\label{prop:decay}
Let the data be informative for approximate $(\varepsilon_q)_{q\in\mathcal{Q}}$-feedback rectification with poles $(\Lambda_q)_{q\in\mathcal{Q}}$, let $V$ and $F_q$ be as in Definition~\ref{def:rectification}, and set $\btri:=\max_q(\|\Lambda_q\|+\varepsilon_q)$. Then $P=V^{-\top}V^{-1}$ satisfies
$(A+BF_q)^\top P(A+BF_q)\preceq\btri^2P$
for all $(A,B)\in\Sigma_q^{\delta_q}$ and all $q\in\mathcal{Q}$. For arbitrary switching signals and $d_k\in\mathbb R^n$, the solutions of \eqref{eq:with_noise} satisfy
\begin{equation}\label{eq:disturbance}
 \|V^{-1}x_k\|\leq\btri^k\|V^{-1}x_0\|
 +\sum_{j=0}^{k-1}\btri^{k-1-j}\|V^{-1}d_j\|.
\end{equation}
Moreover, $\operatorname{dist}(\lambda,\operatorname{spec}\Lambda_q)\leq\varepsilon_q$ for all $\lambda\in\operatorname{spec}(A+BF_q)$, $(A,B)\in\Sigma_q^{\delta_q}$.

In particular, if $\btri<1$, then $x\mapsto x^\top Px$ is a CQLF for all $A+BF_q$ with $(A,B)\in\Sigma_q^{\delta_q}$, $q\in\mathcal Q$.
\end{proposition}
\begin{proof}
The triangle inequality implies $\|V^{-1}(A+BF_q)V\|\leq\btri$ for all $q\in\mathcal Q$. In particular, this holds for the true $(A_q,B_q)\in\Sigma_q^{\delta_q}$. Squaring the induced-norm estimate yields the Lyapunov inequality. Iterating
$\|V^{-1}x_{k+1}\|\leq\btri\|V^{-1}x_k\|+\|V^{-1}d_k\|$
proves \eqref{eq:disturbance}. If $(\Lambda_q+E)z=\lambda z$, $\|z\|=1$, and $\|E\|\leq\varepsilon_q$, normality gives
$\operatorname{dist}(\lambda,\operatorname{spec}\Lambda_q)\leq\|(\lambda I_n-\Lambda_q)z\|=\|Ez\|\leq\varepsilon_q$.
\end{proof}

\noindent Note that $\btri<1$ means that each $\varepsilon_q$ is smaller than the stability margin $1-\|\Lambda_q\|$.

Recall that a system is called \emph{uniformly practically asymptotically
stable} if its solutions approach a neighbourhood of
the origin uniformly with respect to the admissible
disturbances and switching signals, and this
neighbourhood is stable.

\begin{corollary}[Practical stability]\label{cor:practical}
Let the assumptions of Proposition~\ref{prop:decay} hold with $\btri<1$, and let $\|V^{-1}d_k\|\leq\bar d$ for all $k\in\mathbb N_0$. Then, for all $(A,B)\in\Sigma_q^{\delta_q}$, $q\in\mathcal Q$ and arbitrary switching, the set
$\mathcal D:=\{x\in\R^n:\|V^{-1}x\|\leq\bar d/(1-\btri)\}$
is forward invariant and uniformly exponentially attractive, i.e., \eqref{eq:with_noise} is uniformly practically asymptotically stable. For $d_k=0$, the origin is uniformly globally exponentially stable.
\end{corollary}
\begin{proof}
By \eqref{eq:disturbance},
$\|V^{-1}x_k\|\leq\btri^k\|V^{-1}x_0\|+\tfrac{1-\btri^k}{1-\btri}\bar d$,
which yields all claims.
\end{proof}

\begin{remark}\label{rem:beta-dir}
The bound $\btri$ in Proposition~\ref{prop:decay} can be conservative.
A possibly smaller contraction factor is given by
\[
\beta_{\mathrm{dir}}
:=\max_{q\in\Q}\sup_{(A,B)\in\Sigma_q^{\delta_q}}
\|V^{-1}(A+BF_q)V\|
\leq\btri.
\]
For fixed $V$ and feedback gains $F_q$, let $H_q$
be as in Theorem~\ref{thm:lmi-rectification} and set
$N_q=V^{-1}X_q^+K_qH_q$.
By Proposition~\ref{thm:informativity},
\[
\beta_{\mathrm{dir}}
=\max_{q\in\Q}\max_{\|\Delta\|\leq1}
\|N_q-V^{-1}\Omega_q^{1/2}\Delta H_q\|,
\]
where $\Delta\in\R^{n\times r_q}$.
If $\Omega_q\neq0$ for all $q\in\Q$,
Lemma~\ref{lem:norm} shows that
$\beta_{\mathrm{dir}}=\sqrt{\gamma^\star}$,
where $\gamma^\star$ is the optimal value of
\[
\begin{aligned}
\min_{\gamma,t_1,\ldots,t_Q}\quad &\gamma\\
\text{subject to}\quad
&\mathcal{M}_q(\gamma,N_q,V,H_q,t_q)\succeq0,\,t_q>0,\, q\in\Q.
\end{aligned}
\]
Since $V$ and the gains $F_q$ are fixed,
this is a semidefinite optimization problem
with only scalar decision variables.
In particular, $\beta_{\mathrm{dir}}<1$
is sufficient for all conclusions of
Proposition~\ref{prop:decay} and Corollary~\ref{cor:practical}, with
$\btri$ replaced by $\beta_{\mathrm{dir}}$.
Unlike $\btri$, $\beta_{\mathrm{dir}}$ measures how much of the stability margin $1-\|\Lambda_q\|$ is actually consumed by the uncertainty.
\end{remark}

\section{Numerical example}\label{sec:examples}
\begin{figure}[tbp]
\centering
\includegraphics[width=\linewidth]{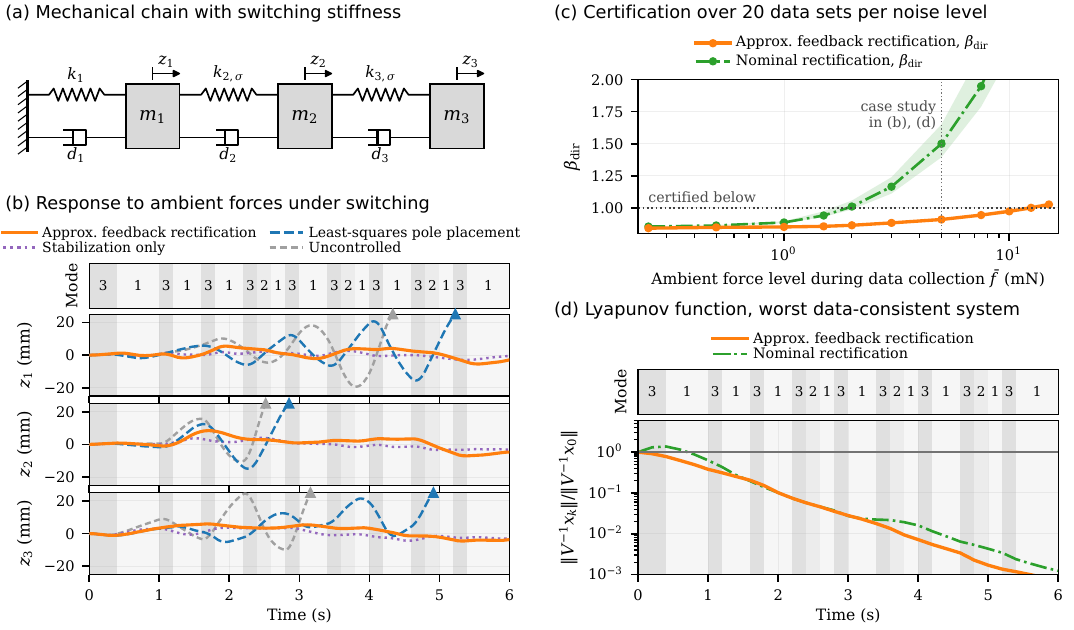}
\caption{(a) Chain with switching stiffnesses. (b) Displacements of the true system from rest under ambient forces during operation, drawn uniformly from $[-50,50]$\,mN on each mass, compared to $\bar f=5$\,mN during data collection. The triangle markers indicate where a response leaves the plotted range. (c) $\bdir$ over the ambient force level during data collection, median and range over $20$ data sets. The certified force level is more than six times larger with approximate feedback rectification. (d) $\|V^{-1}x_k\|$ on the data-consistent system that maximizes the one-step growth, worst initial state, no ambient forces.}
\label{fig:num}
\end{figure}

\emph{Setup.} We consider Example~2 of \cite{WeiMottershead16}, which is a mass-spring-damper chain, a common building block of mechatronic systems, that consists  of three masses $m_1$, $m_2$ and $m_3$. Here $m_1$ is attached to a wall and $m_3$ is free, see Fig.~\ref{fig:num}(a). In addition, actuator forces $\mathcal Bu$ and ambient forces $f$ act along the chain, which gives one system per mode
\[
\ddot z+c\,\mathcal D\dot z+\mathcal K_{q} z=\mathcal Bu+f,\qquad q\in\Q=\{1,2,3\},
\]
with $\mathcal D$, $\mathcal B$, and the structure of $\mathcal K_q$ as in \cite{WeiMottershead16}, $k_1=5$\,N/m, $(k_{2,q})_q=(5,6.5,9.86)$\,N/m, and $(k_{3,q})_q=(5,12.5,29.3)$\,N/m. Mode~1 is the original system, $c=0.2$ makes the chain lightly damped, $f$ is an unmeasured disturbance. The stiffness switches only at the sampling instants, so zero-order-hold sampling with $\tau=0.2$\,s of $x=(z,\dot z)$ in cm and cm/s yields $x_{k+1}=A_qx_k+B_qu_k+E_qf_k$ in mode $q$ with state and input dimensions $n=6$ and $m=3$. Here, $E_qf_k$ is the disturbance. Each mode is stable, but switching between soft mode~1 and stiff mode~3 can inject energy into the chain and destabilize it.

\emph{Specification.} In every mode, the closed loop shall have the natural frequencies $\omega\in\{1.00,2.79,4.03\}$\,rad/s of mode~1 and the decay rate $1.5\,\mathrm s^{-1}$, i.e., $\Lambda_q=\Lambda$ with eigenvalues $e^{\tau(-1.5\pm\mathrm i\omega)}$, that is $\|\Lambda\|=0.741$.

\emph{Data.} The chain is operated with the velocity feedback $u=K_0x$, $u=-0.05\,\mathcal B^\top\dot z$ with $\dot z$ in m/s. In this closed loop, only actuators~1 and~2 are additionally excited, $u_k=K_0x_k+J\xi_k$ with the unit vectors $J=[e_1,e_2]$ and $\xi_k$ uniformly distributed on $[-1,1]^2$\,N. In every mode~$q$, $T_q=120$ samples of $x_k$, $u_k$, and $x_{k+1}$ form $X_q$, $U_q$, and $X_q^+$. If the ambient forces satisfy $(\sum_k\|f_k\|^2/T_q)^{1/2}\le\bar f$, then $\|W_q\|\le\delta_q:=\bar e\sqrt{T_q}\,\bar f$ with $\bar e:=21\ge\|E_q\|$. The case study uses $\bar f=5$\,mN; Fig.~\ref{fig:num}(c) varies~$\bar f$.

\emph{Designs.} Since $U_q=K_0X_q+J\Xi_q$, where $\Xi_q$ collects the $\xi_k$, we have $\rank\smat{X_q\\ U_q}=n+2=8<n+m$, and $(A_q,B_q)$ cannot be uniquely identified. By Proposition~\ref{thm:informativity}, every gain with a finite worst-case error is of the form $F_q=U_qK_qH_qV^{-1}$ with $X_qK_qH_q=V$, hence $F_q=K_0+JG_q$ with $G_q=\Xi_qK_qH_qV^{-1}$: the data only allow changes in the excited input directions. The rectification designs have this form; the baselines are designed on the least-squares fit of $(A_q+B_qK_0,B_qJ)$ and thus also retune only actuators~1 and~2. \emph{Least-squares pole placement} assigns $\spec\Lambda$ to this fit in every mode. \emph{Stabilization only} computes, for the same fits, a common $P_{\mathrm s}\succ0$ with $(A_q+B_qF_q)^\top P_{\mathrm s}(A_q+B_qF_q)\preceq\|\Lambda\|^2P_{\mathrm s}$. \emph{Approximate feedback rectification} minimizes $\varepsilon=\varepsilon_q$ subject to Corollary~\ref{cor:lmi-reference}, with $V_0$ from the pole-placement eigenvectors of mode~1, and re-solves with $V_0$ replaced by $V$ up to $20$ times. \emph{Nominal rectification} does the same with $\Omega_q=0$. 

\begin{table}[t]
\centering
\caption{Closed loop systems. $\Delta\lambda$: largest eigenvalue error.}
\label{tab:num}
\setlength{\tabcolsep}{4pt}
\begin{tabular}{@{}lccc@{}}
\toprule
 & $\Delta\lambda$ & certificate & $\operatorname{cond.} V$\\
\midrule
Least-squares pole placement & $0.00$ & -- & --\\
Stabilization only & $0.70$  & common $P_{\mathrm s}$ & --\\
Nominal rectification & $0.06$ & $\bdir=1.58$ & $273$\\
Approx.\ feedback rectification & $0.11$ & $\bdir=0.92$ & $40$\\
\bottomrule
\end{tabular}
\end{table}

\emph{Results.} Least-squares pole placement meets the specification (Table~\ref{tab:num}), but the switching signal in Fig.~\ref{fig:num}(b) destabilizes it, as well as the uncontrolled chain. Stabilization only is certified but misses the prescribed poles; its highest frequency varies between $5.7$ and $8.0$\,rad/s. Approximate feedback rectification achieves $\varepsilon=0.185$, hence $\btri=\|\Lambda\|+\varepsilon=0.925<1$, and Proposition~\ref{prop:decay} yields a CQLF for all data-consistent systems, eigenvalues within $\varepsilon$ of $\spec\Lambda$, and, with Corollary~\ref{cor:practical}, bounded responses to ambient forces. Nominal rectification is closer to the specification on the true system, but its ill-conditioned $V$ amplifies the uncertainty, $\bdir=1.58$. In Fig.~\ref{fig:num}(d), $\|V^{-1}x_k\|$ overshoots its initial value for a data-consistent system, so $x^\top V^{-\top}V^{-1}x$ is no valid Lyapunov certificate; for approximate feedback rectification it decreases monotonically. Fig.~\ref{fig:num}(c): certified up to $\bar f=1.5$\,mN (nominal) and $10$\,mN (approximate).

\section{Conclusion}
We characterized when disturbed and possibly rank-deficient input-state data are informative for approximate feedback rectification and derived an LMI synthesis of the common eigenvectors and gains. If the rectification error is smaller than the stability margin of the prescribed poles, the data certify stability under arbitrary switching.

\section*{AI Disclosure}
The authors used AI tools (ChatGPT Astra Pro, Claude Opus 5.5) for verifying the mathematical results and improving the writing. The code for the simulations in Section~\ref{sec:examples} was developed with the assistance of these tools. The authors checked all AI-generated content and take full responsibility for it.

\phantomsection
\pdfbookmark[1]{References}{archive-references}
\bibliographystyle{IEEEtran}
\bibliography{references}
\end{document}